\documentclass[12pt]{amsart}
\usepackage{amsmath,amssymb,latexsym}

\theoremstyle{plain}
\newtheorem{Thm}{Theorem}
\newtheorem{Cor}{Corollary}
\newtheorem{Prop}{Proposition}[section]
\newtheorem{Lem}[Prop]{Lemma}
\newtheorem{Def}{Definition}

\input{epsf}
\theoremstyle{definition}

\begin{document}

\title
{The Kauffman polynomials of generalized Hopf links}

\author{Jianyuan K. Zhong}
\address{Department of Mathematics \& Statistics\\
   Louisiana Tech University\\
    Ruston, LA 71272}
\email{kzhong@coes.latech.edu}
\author{Bin Lu}
\address{Department of Mathematics \\
   The University of Arizona\\
    Tucson, AZ 85721} 
\email{binlu@math.arizona.edu} 
\date{November 28, 2001}
\begin{abstract}
Following the recent work by Chan \cite{C00} and Morton and Hadji
\cite{MH01} on  the Homflypt polynomials of some generalized Hopf
links, we investigate the  Kauffman polynomials
of generalized Hopf links. By studying the Kauffman skein
module of the solid torus $S^1\times D^2$, we establish a
similar skein map on the Kauffman skein module of $S^1\times D^2$
which has distinct eigenvalues. Furthermore we are able to
calculate the Kauffman polynomials of some specific generalized Hopf links.
\end{abstract}

\maketitle
\section{Introduction}

In \cite{C00}, Chan discusses the Homflypt polynomials of reversed
string parallels $H(k_1, k_2; n_1, n_2)$ of the Hopf link. Morton
and Hadji analyze their structure more closely by using the Homflypt
skein of the annulus and identify the eigenvalues and eigenvectors
(see \cite{MH01}). Here we
 follow a similar approach to
Morton and Hadji's to investigate Kauffman polynomials of
generalized Hopf links by using the Kauffman skein module  of the solid torus (the annulus). 

Let $k$ be an integral domain containing the invertible elements
$\alpha $ and $s$. We assume that $s-s^{-1}$ is invertible in $k$.

By a {\it framed oriented link} we mean a link equipped with a
string orientation together with a nonzero normal vector field up
to homotopy. By a {\it  framed link} we mean an unoriented framed
link. The links described by figures in this paper will be
assigned the vertical framing which points towards the reader.

There are various versions of {\it the Kauffman polynomial} in
the literature. The version we use  in this paper is defined by
the following Kauffman skein relations:
$$
\raisebox{-3mm}{\epsfxsize.3in\epsffile{left1.ai}}\quad
        -\quad\raisebox{-3mm}{\epsfxsize.3in\epsffile{right1.ai}}\quad
=\quad (\ s - \
s^{-1})\bigl(\quad\raisebox{-3mm}{\epsfxsize.3in\epsffile{parra1.ai}}-\quad\raisebox{-3mm}{\epsfxsize.3in\epsffile{parra2.ai}}\quad\bigr)\quad
,
$$
$$
\raisebox{-3mm}{\epsfxsize.35in\epsffile{framel1.ai}}
        \quad=\quad \alpha \quad \raisebox{-3mm}{\epsfysize.3in\epsffile{orline1.ai}}\quad  ,
$$
$$
L\ \sqcup
\raisebox{-2mm}{\epsfysize.3in\epsffile{unknot1.ai}}\quad =\quad
\delta\quad L\quad.
$$
where $\delta=({\dfrac{\alpha-\alpha^{-1}} {\ s - \ s^{-1}}}+1)$. The last relation follows from the first two when $L$ is nonempty.

\noindent {\bf Remark.}
The Kauffman polynomial of the empty link is normalized  to be $1$ and
the Kauffman polynomial of a link $L$ is denoted by $<L>$.

Let $M$ be a smooth, compact and oriented $3$-manifold.
\begin{Def}
{\bf The Kauffman skein module} of $M$, denoted by $K(M)$, is the
$k$-module freely generated by isotopy classes of framed links in
$M$ including the empty link quotient by the Kauffman skein
relations given above.
\end{Def}

We consider the unoriented Hopf link with linking number $1$ as
shown below:
$$\hbox{Hopf link}=\raisebox{-7mm}{\epsfysize0in\epsffile{hopf.ai}}.$$
We have that $<\hbox{Hopf link}>= \delta[\delta+(s-s^{-1})(\alpha-\alpha^{-1})]$.

A great variety of links can be realized through the construction
of satellite links by using the Hopf link \cite{MH01}, which is to
take the two components of the Hopf link and decorate them with
diagrams in the annulus. If the two components are decorated by
$P_1$ and $P_2$, respectively, we denote the realized satellite
link by $H(P_1, P_2)$. When $P_1$ is $n$ parallel copies of the
unknot and $P_2$ is $m$ parallel copies of the unknot, we denote
$H(P_1, P_2)$ by $H(n, m)$. These are the generalized Hopf links.
We will label a string with an integer $n$ to indicate $n$
parallel strings.

A picture of $H(n, m)$ is given by the following:
$$H(n,m)=\raisebox{-18mm}{\epsfysize0in\epsffile{2nk.ai}}\quad=\quad\raisebox{-21mm}{\epsfysize0in\epsffile{1nk.ai}}\quad .$$We observe that the links $H(n, m)$ and $H(m, n)$ are equivalent
links.
\section{On the Kauffman skein module of $S^1\times D^2$}

\subsection{Hecke algebras}
Here we summarize a geometric realization of the Hecke algebras through the Homflypt skein modules.
\begin{Def}
{\bf The Homflypt skein module} of $M$, denoted by $S(M)$,
is the $k$-module freely generated by isotopy classes of framed
oriented links in $M$ including the empty link, quotient  by the
Homflypt skein relations given in the following figure.
$$
x^{-1}\quad\raisebox{-3mm}{\epsfxsize.3in\epsffile{left.ai}}\quad
        -\quad x\quad\raisebox{-3mm}{\epsfxsize.3in\epsffile{right.ai}}\quad
=\quad (\ s - \
s^{-1})\quad\raisebox{-3mm}{\epsfxsize.3in\epsffile{parra.ai}}\quad
,
$$
$$
\raisebox{-3mm}{\epsfxsize.35in\epsffile{framel.ai}}
        \quad=\quad (xv^{-1})\quad \raisebox{-3mm}{\epsfysize.3in\epsffile{orline.ai}}\quad  ,
$$
$$
L\ \sqcup \raisebox{-2mm}{\epsfysize.3in\epsffile{unknot.ai}}\quad
                                        =\quad {\dfrac{v^{-1}-v} {\ s - \ s^{-1}}}\quad L\quad.
$$
\end{Def}

\begin{Def}

{\bf The relative Homflypt skein module}. Let $X=\{x_1,
x_2,\cdots, x_{n}\}$ be a finite set of framed points oriented
negatively (called input points) in the boundary $\partial M$, and let
$Y=\{y_1, y_2,\cdots, y_{n}\}$ be a finite set of framed points
oriented positively (called output points) in
$\partial M$. The relative skein module $S(M, X, Y)$ is
the $k$-module generated by relative framed oriented links in
$(M,\partial M)$ such that $L \cap \partial M=\partial L=\{x_i,
y_i\}$ with the induced framing and orientation, considered up to
an ambient isotopy fixing $\partial M$, and quotient  by the Homflypt
skein relations.

\end{Def}

In the cylinder $D^2\times I$, let $X_n$ be a set of $n$ distinct
input framed points on a diameter $D^2\times \{1\}$ and $Y_n$ be a
set of $n$ distinct output framed points on a diameter $D^2\times
\{0\}$, it is a well-known result \cite{AM98} that the relative
Homflypt skein module $K(D^2\times I, X_n\amalg Y_n)$ is
isomorphic to the {\it $n$th Hecke algebra} $H_n$, which is the
quotient of the braid group algebra $k[B_n]$ by the Homflypt skein
relations. This is the geometric realization of $H_n$.

In $H_n$, there is a set of quasi-idempotent elements $e_\lambda$s
\cite{G86} \cite{AM98} \cite{cB98}. After normalization, we denote
the corresponding idempotent by $y_\lambda$ \cite{cB98}.

\subsection{The Birman-Murakami-Wenzl category and the
Birman-Murakami-Wenzl algebra $K_n$}

In this section 2.2 and section 2.3, we give a summary of the 
related work of  Beliakova
and Blanchet. Details, further references to the origin of some of
these ideas, and related results of others can be found in
\cite{BB2000}. We also  provide some figures to illustrate the ideas.

\begin{Def}
{\bf The relative Kauffman skein module}. Let $X=\{x_1,
x_2,\cdots, x_{n}, y_1, y_2,\cdots, y_{n}\}$ be a finite set of $2n$
framed points in the boundary $\partial M$.  The relative
skein module $K(M, X)$ is the $k$-module generated by relative
framed links in $(M,\partial M)$ such that $L \cap
\partial M=\partial L=\{x_i, y_i\}$ with the induced framing,
considered up to an ambient isotopy fixing $\partial M$,  quotient  by the
Kauffman skein relations.
\end{Def}

{\it The Birman-Murakami-Wenzl category $K$} consists of: an
object in $K$ is a disc $D^2$ equipped with a finite set of points
and a nonzero vector at each point. If $\beta=(D^2, l_0)$ and
$\gamma=(D^2, l_1)$ are objects, the module $Hom_K(\beta, \gamma)$
is $K(D^2\times [0,1], l_0\times 0 \amalg l_1\times 1)$. For a
Young diagram $\lambda$, we denote by $\square_\lambda$ the object
of the category $K$ formed with one point assigned for
 each cell of $\lambda$. We will
use the notation $K(\beta, \gamma)$ for $Hom_K(\beta, \gamma)$.
For composition of morphisms $f$ and $g$, it's by stacking $f$ on the top of $g$:
$$K(\beta, \gamma)\times K(\gamma, \delta)\to K(\beta, \delta),$$
$$(f,g)\to fg.$$
Note that Beliakova and Blanchet choose to stack the second one on
the top of the first. Here we follow the convention as in \cite{ZG99} \cite{AM98}.

As a special case, in the cylinder $D^2\times I$, let $X_n$ be a
set of $n$ distinct framed points on a diameter $D^2\times \{1\}$
and $Y_n$ be a set of $n$ distinct framed points on a diameter
$D^2\times \{0\}$, then the relative Kauffman skein module
$K(D^2\times I, X_n\amalg Y_n)$ is isomorphic to the {\it
Birman-Murakami-Wenzl algebra} $K_n$, which is the quotient of the
braid group algebra $k[B_n]$ by the Kauffman skein relations.

The Birman-Murakami-Wenzl algebra $K_n$ is
generated by the identity ${\mathbf{1}}_n$, positive
transpositions $e_1, e_2, \cdots, e_{n-1}$ and hooks $h_1, h_2,
\cdots,$ $h_{n-1}$ as the following:
$$e_i=\quad \raisebox{-3.5mm}{\epsfxsize.0in\epsffile{ei.ai}}$$
$$h_i=\quad \raisebox{-3.5mm}{\epsfxsize.0in\epsffile{hi.ai}}$$
for $1\leq i\leq n-1$.

Then $K_n$ is the braid group algebra $k[B_n]$ quotient  by the
following relations:

\hspace{2cm} $(B_1)\ e_ie_{i+1}e_i=e_{i+1}e_{i}e_{i+1}$,

\hspace{2cm} $(B_2)\ e_ie_j=e_je_{i},\ |i-j|\geq 2$,

\hspace{2cm} $(R_1)\ h_i e_i=\alpha^{-1}h_i$,

\hspace{2cm} $(R_2)\ h_i{e_{i-1}}^{\pm 1}h_i={\alpha}^{\pm 1}h_i$,

\hspace{2cm} $(K)\ e_i-e_i^{-1}=(s-s^{-1})({\mathbf 1}_n-h_i)$.

Let $I_n$ be the ideal generated by $h_{n-1}$. Note that
$I_n=\{(a\otimes {\mathbf 1}_1)h_{n-1}(b\otimes {\mathbf 1}_1): a,
b\in K_{n-1}\}$. The quotient of $K_n$ by the ideal $I_n$ is
isomorphic to the $n$th Hecke algebra $H_n$. We denote the canonical projection map by
$\pi_n$:
$$\pi_n: K_n\to H_n.$$

\begin{Thm}
Over ${\mathbb Q}(\alpha, s)$, the field
of rational functions in $\alpha, s$, there exists a multiplicative
homomorphism $s_n:H_n\to K_n$, such
that
$$\pi_n\circ s_n=id_{H_n},$$
$$s_n(x)y=ys_n(x)=0, \ \forall x\in H_n, \forall y \in I_n.$$
\end{Thm}

\begin{Cor}
$K_n\cong H_n\oplus I_n$.
\end{Cor}

Let $\widetilde{y}_\lambda=s_n(y_\lambda)$.

\subsection{A basis for the Birman-Murakami-Wenzl algebra
$K_n$}

We call a sequence $\Lambda=(\Lambda_1, \cdots, \Lambda_n)$ of Young
diagrams an up and down tableau of length $n$ and
shape $\Lambda_n$ if two consecutive Young diagrams $\Lambda_i$
and $\Lambda_{i+1}$ differ by exactly one cell. We observe that in an up
and down tableau $\Lambda=(\Lambda_1, \cdots, \Lambda_n)$ of
length $n$, the size of $\Lambda_n$ is either $n$ or less than $n$
by an even number.

For an up and down tableau $\Lambda$ of length $n$, we denote by
$\Lambda'$ the tableau of length $n-1$ obtained by removing the
last Young diagram in the sequence $\Lambda$. We define
$a_{\Lambda}\in K(n, \square_{\Lambda})$ and $b_{\Lambda}\in K(
\square_{\Lambda}, n)$ by
$$a_1=b_1={\mathbf{1}}_1.$$
If $|\Lambda_n|=|\Lambda_{n-1}|+1$, then
$$a_{\Lambda}=(a_{\Lambda'}\otimes
{\mathbf{1}}_1){\widetilde{y}_{\Lambda_n}},$$
$$b_{\Lambda}={\widetilde{y}_{\Lambda_n}}(b_{\Lambda'}\otimes
{\mathbf{1}}_1);$$ if $|\Lambda_n|=|\Lambda_{n-1}|-1$, then
$$a_{\Lambda}=\frac{<\Lambda_n>}{<\Lambda_{n-1}>}(a_{\Lambda'}\otimes
{\mathbf{1}}_1)({\widetilde{y}_{\Lambda_n}}\otimes \cup),$$
\hspace{3.5cm} $b_{\Lambda}=({\widetilde{y}_{\Lambda_n}}\otimes
\cap)(b_{\Lambda'}\otimes {\mathbf{1}}_1).$

Here $<\lambda>$ is the quantum dimension \cite{W88} associated with
$\lambda$, which is the Kauffman polynomial of
$\widehat{\widetilde{y}}_\lambda$ in $S^3$. Note that $<\lambda>$ is
invertible in ${\mathbb Q}(\alpha, s)$ \cite{BB2000}[Theorem 7.5].

We provide the following figures to illustrate the ideas.

(1) If $|\Lambda_n|=|\Lambda_{n-1}|+1$, then
$$a_{\Lambda}=\raisebox{-20mm}{\epsfxsize.0in\epsffile{ala.ai}}\quad,\quad  \quad
b_{\Lambda}=\raisebox{-20mm}{\epsfxsize.0in\epsffile{bla.ai}}.$$

(2) If $|\Lambda_n|=|\Lambda_{n-1}|-1$, then
$$a_{\Lambda}=\frac{<\Lambda_n>}{<\Lambda_{n-1}>}\quad
\raisebox{-20mm}{\epsfxsize.0in\epsffile{ala1.ai}}\quad, \quad
b_{\Lambda}=\raisebox{-20mm}{\epsfxsize.0in\epsffile{bla1.ai}}.$$

\begin{Thm}
The family $a_{\Lambda}b_{\Gamma}$ for all up and down tableaux
$\Lambda, \Gamma$ of length $n$ such that $\Lambda_n=\Gamma_n$
forms a basis for $K_n$.
\end{Thm}

(1) If $|\Lambda_n|=|\Lambda_{n-1}|+1$ (so is
$|\Gamma_n|=|\Gamma_{n-1}|+1$), then
$$a_{\Lambda}b_{\Gamma}=\raisebox{-30mm}{\epsfxsize.0in\epsffile{alabth.ai}}.$$
(2) If $|\Lambda_n|=|\Lambda_{n-1}|-1$ (so is
$|\Gamma_n|=|\Gamma_{n-1}|-1$), then
$$a_{\Lambda}b_{\Gamma}=\frac{<\Lambda_n>}{<\Lambda_{n-1}>}\quad
\raisebox{-35mm}{\epsfxsize.0in\epsffile{alabth1.ai}}.$$

Let $\Lambda=(\Lambda_1, \cdots, \Lambda_n), \
\Gamma=(\Gamma_1,\cdots, \Gamma_n)$ be two up and down tableaux of
length $n$.
If $\Lambda=\Gamma$, i.e. $\Lambda_i=\Gamma_i$ for $1\leq i \leq
n$, then $b_{\Gamma}a_{\Lambda}=\widetilde{y}_{\Lambda_n}$;
otherwise $b_{\Lambda}a_{\Gamma}=0$. This follows from the
corresponding properties in the Hecke category \cite{cB98}. These
properties will be used in the following section.

\subsection{The Kauffman skein module of the solid torus $S^1\times
D^2$}

There is a natural wiring of the cylinder $D^2\times I$ into the
solid torus $S^1\times D^2$,
$$\quad \raisebox{-15mm}{\epsfxsize.0in
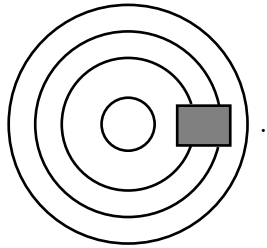
\epsffile{wiring.ai}}\ .$$
We denote the image of $K_n$ under the above wiring by $\widehat{K_n}$.

\begin{Thm}
$\widehat{K_n}\subseteq K(S^1\times D^2)$ has a basis given by the
collection \{$\widehat{\widetilde{y}}_{\lambda}: |\lambda|$ is either $n$
 or less than $n$ by an even number\}.
\end {Thm}

\begin{proof}
By Theorem 2 in the previous section, $K_n$ has a basis given by
the family $a_\Lambda b_\Gamma$, so $\widehat{K_n}$ is generated
by the family $\widehat{a_\Lambda b_\Gamma}$. Since
$\widehat{a_\Lambda b_\Gamma}=\widehat{b_\Gamma
a_\Lambda}=\delta_{\Lambda\Gamma}\widehat{b_\Lambda
a_\Lambda}=\delta_{\Lambda\Gamma}\widehat{\widetilde{y}}_{\Lambda_n}$,
where $\Lambda_n$ is a Young diagram of size either $n$ or less
than $n$ by an even number. Thus $\widehat{K_n}$ is generated by
the collection \{$\widehat{\widetilde{y}}_{\lambda}: |\lambda|$ is
either $n$ or less than $n$ by an even number \}.

Now we want to prove the above generating set is linearly
independent. We show this by comparing the dimensions.

From Corollary 1, we have $K_n\cong H_n\oplus I_n$, so
$\widehat{K_n}\cong \widehat{H_n}+\widehat{I_n}$. As
$I_n=\{(a\otimes {\mathbf 1}_1)h_{n-1}(b\otimes {\mathbf 1}_1):\ a,
b\in K_{n-1}\}$, a typical element in $\widehat{I_n}$ looks like
the following:
$$\quad \raisebox{-3.5mm}{\epsfxsize.0in
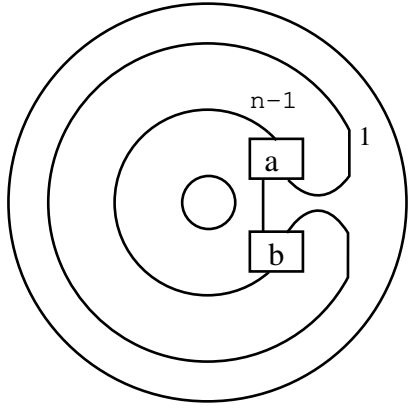
\epsffile{hatin.ai}}$$
where $a, b \in K_{n-1}$.
We can see that the above element is in $\widehat{K_{n-2}}$. On the other hand, we
have $\widehat{K_{n-2}}\subseteq \widehat{I_n}$, therefore
$\widehat{I_n}=\widehat{K_{n-2}}$. i.e. $\widehat{K_n}\cong
\widehat{H_n}+\widehat{K_{n-2}}$.

Repeating the process for $\widehat{K_{n-2}}$, we conclude that

\begin{equation}
    \widehat{K_n}\cong
    \begin{cases}
       \widehat{H_1}+\widehat{H_3}+\cdots+\widehat{H_n}, &\textrm {if $n$ is odd;}\notag\\
       <\phi>+\widehat{H_2}+\widehat{H_4}+\cdots+\widehat{H_n},  &\textrm{if $n$ is
       even;}\notag
    \end{cases}
\end{equation}
where $\phi$ is the empty link.

Since $\widehat{H_i}\cap \widehat{H_j}=0$ whenever $i\neq j$ in
the Homflypt skein module of the solid torus, the above
decomposition is a direct sum.

\begin{equation}
    \widehat{K_n}\cong
    \begin{cases}
       \widehat{H_1}\oplus\widehat{H_3}\oplus\cdots\oplus\widehat{H_n},  &\textrm {if $n$ is odd;}\notag\\
       <\phi>\oplus\widehat{H_2}\oplus\widehat{H_4}\oplus\cdots\oplus\widehat{H_n}, &\textrm{if $n$ is
       even.}\notag
    \end{cases}
\end{equation}
Therefore we have the following equality for the dimensions.

\begin{equation}
    \dim(\widehat{K_n})=
    \begin{cases}
       \dim(\widehat{H_1})+\dim(\widehat{H_3})+\cdots+\dim(\widehat{H_n}),   &\textrm {if $n$ is odd;}\notag\\
       1+\dim(\widehat{H_2})+\dim(\widehat{H_4})+\cdots+\dim(\widehat{H_n}),  &\textrm{if $n$ is even.}
    \end{cases}
\end{equation}

As we know from \cite{AM98} that the dimension of
$\widehat{H_i}$($=C_i$) in $S(S^1\times D^2)$ is equal to the number
of Young diagrams of size $i$. By comparing the number of generators
in the set \{$\widehat{\widetilde{y}}_{\lambda}: |\lambda|$ is either
$n$ or less than $n$ by an even number \} and the
dimension of $\widehat{K_n}$, they are equal, so the given generating set must be linearly independent. Thus it forms a basis for
$\widehat{K_n}$.

\end{proof}

\begin{Cor}
The collection of all the elements
\{$\widehat{\widetilde{y}_{\lambda}}: \lambda$ is any Young
diagram\} forms a basis for $K(S^1\times D^2)$ over ${\mathbb
Q}(\alpha, s)$,
where
$$\widehat{\widetilde{y}}_\lambda=\raisebox{-24mm}{\
\epsfxsize1.9in
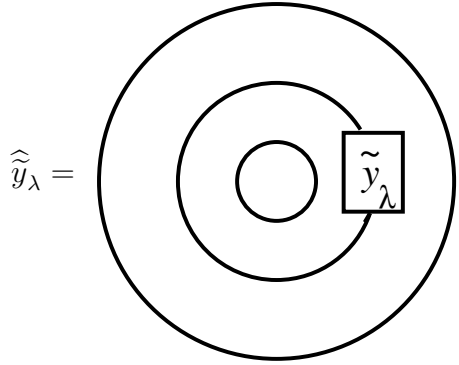
\epsffile{yla1.ai}}.$$
\end{Cor}
\begin{proof}
Let $L$ be a framed link in $S^1\times D^2$, up to a scalar
multiple, then $L$ is the closure of an $n$-strand braid for some
integer $n\geq 0$ \cite{PS}[6.5 Alexander's Braiding Theorem]; the
braid quotient by the Kauffman skein relations is in $K_n$.
Therefore $L\in \widehat{K_n}$ in $K(S^1\times D^2)$. This shows
any element in $K(S^1\times D^2)$ lies in $\bigcup_{n\geq
0}\widehat{K_n}$, i.e. $K(S^1\times D^2)\subseteq\bigcup_{n\geq
0}\widehat{K_n}$, hence $K(S^1\times D^2)=\bigcup_{n\geq
0}\widehat{K_n}$. So we need only to show that the set
\{$\widehat{\widetilde{y}}_{\lambda}: \lambda$ is a Young
diagram\} forms a basis for $\bigcup_{n\geq 0}\widehat{K_n}$.

From the previous Theorem, $\widehat{K_n}\subseteq K(S^1\times
D^2)$  has a basis given by the collection
\{$\widehat{\widetilde{y}}_{\lambda}: |\lambda|$ is either $n$ or
less than $n$ by an even number \}. It follows that the set
\{$\widehat{\widetilde{y}}_{\lambda}: 0\leq |\lambda|\leq n$\} forms a basis
for $\bigcup_{0\leq i\leq n}\widehat{K_i}$. As we have
$K(S^1\times D^2)=\bigcup_{n\geq 0}\widehat{K_n}$. The corollary
follows by induction on $m$ for $\bigcup_{0\leq n\leq
m}\widehat{K_n}$.

\end{proof}

\subsection{A linear map on $K(S^1\times D^2)$}

We define a natural linear map $\varphi$ on $K(S^1\times D^2)$ in
the following:
$$\varphi: K(S^1\times D^2)\to K(S^1\times D^2)$$
$$\raisebox{-9.5mm}{\epsfxsize.0in
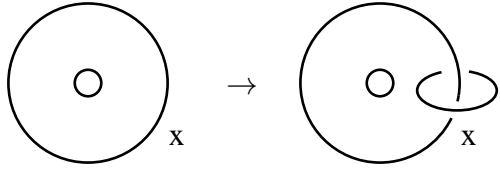
\epsffile{3map.ai}}\quad\to\quad
\raisebox{-9.5mm}{\epsfxsize.0in\epsffile{4map.ai}}.$$
When
$X=A_1$, where $A_1$ is the unknot in the annulus representing the
longitude, then $H(n, m)=\varphi^m(A_1)$.

Let $\lambda$ be a Young diagram, $\lambda$ has $|\lambda|$ cells indexed by \{$(i, j)$\}. If $c$ is the cell of index $(i, j)$ in $\lambda$, its content $cn(c)$ is defined by
$cn(c)=j-i$. We define,
$$c_\lambda=\delta+(\ s - \
s^{-1})(\alpha\sum_{c\in \lambda}s^{2cn(c)}-\alpha^{-1}\sum_{c\in
\lambda}s^{-2cn(c)}).$$

Note that the scalars $c_\lambda$ are different for each Young
diagram $\lambda$.
\begin{Prop}
$$\varphi(\widehat{\widetilde{y}}_\lambda)=c_\lambda \widehat{\widetilde{y}}_\lambda.$$
\end{Prop}
The proof follows from the following lemma.

\begin{Lem}
$$\raisebox{-15mm}{\epsfxsize.0in\epsffile{alasab.ai}}\quad=c_\lambda \quad
\raisebox{-12.5mm}{\epsfxsize.0in\epsffile{alasa.ai}}.$$
\end{Lem}

\begin{proof}
We proceed by induction on $|\lambda|$.

(1) It's easy to check the result is true when $|\lambda|=0$ and $|\lambda|=1$.

(2) Suppose it is true for $|\lambda|<n$. Now assume $\lambda$ is a Young diagram of size $n$.

First we have the following identity by the absorbing property given in \cite{BB2000}[Chapter 5]: $$\widetilde{y}_\lambda=\widetilde{y}_\lambda(\widetilde{y}_{\lambda'}\otimes 1)\widetilde{y}_\lambda,$$
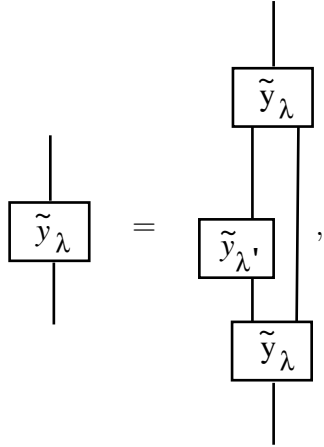
$$\raisebox{-12mm}{\epsfxsize.0in\epsffile{alasa.ai}}\quad=\quad
\raisebox{-28mm}{\epsfxsize.0in\epsffile{alas1.ai}},$$
and the skein relation \cite[Prop. 6.1]{BB2000}:
$$\raisebox{-27mm}{\ \epsfxsize0in\epsffile{alasb.ai}}
=s^{2cn(c)}\quad \raisebox{-12mm}{\
\epsfxsize0in\epsffile{alasa.ai}},$$
and from the above skein relation, we can develop the following relation,
$$\raisebox{-27mm}{\ \epsfxsize0in
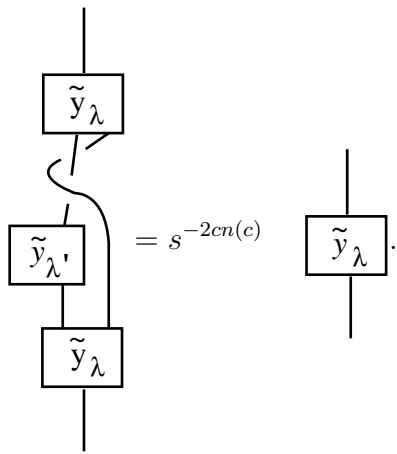
\epsffile{alasb1.ai}}
=s^{-2cn(c)}\quad \raisebox{-12mm}{\
\epsfxsize0in\epsffile{alasa.ai}}.$$
where $cn(c)$ is the content
of the extreme cell $c$ of $\lambda$ to be removed to obtain
$\lambda'$.

Now, apply
the Kauffman skein relation to the following,
$$\raisebox{-28mm}{\epsfxsize.0in\epsffile{alas2.ai}}\quad=\quad
\raisebox{-28mm}{\epsfxsize.0in\epsffile{alas3.ai}}+(s-s^{-1})(
\raisebox{-28mm}{\epsfxsize.0in\epsffile{alas4.ai}}\quad+\quad
\raisebox{-28mm}{\epsfxsize.0in\epsffile{alas5.ai}})$$
$$=\quad
\raisebox{-28mm}{\epsfxsize.0in\epsffile{alas3.ai}}+(s-s^{-1})(\alpha
s^{2cn(c)}-\alpha^{-1}s^{-2cn(c)})
\raisebox{-12mm}{\epsfxsize.0in\epsffile{alasa.ai}}$$
(by induction on $\lambda'$)
$$=c_{\lambda'}\quad
\raisebox{-28mm}{\epsfxsize.0in\epsffile{alas1.ai}}+(s-s^{-1})(\alpha
s^{2cn(c)}-\alpha^{-1}s^{-2cn(c)})
\raisebox{-12mm}{\epsfxsize.0in\epsffile{alasa.ai}}$$
\hspace{1cm} $=c_\lambda \widetilde{y}_\lambda$.
\end{proof}

\begin{Cor}
For $m\in {\mathbf Z}^+$,
$$\varphi^m (\widehat{\widetilde{y}}_\lambda)=(c_\lambda)^m \widehat{\widetilde{y}}_\lambda.$$
\end{Cor}

\begin{Thm}
The eigenvalues of $\varphi |_{\widehat{K_n}}$ are all distinct.
\end{Thm}
\begin{proof}
Let $\lambda$ be a Young diagram of size $i$. From above, we have
$\varphi (\widehat{\widetilde{y}}_\lambda)=(c_\lambda)
\widehat{\widetilde{y}}_\lambda$, so
$\widehat{\widetilde{y}}_\lambda$ is an eigenvector of $\varphi$
with eigenvalue $c_\lambda$. There are $p(i)$ of these
eigenvectors all with distinct eigenvalues, where $p(i)$ is the
number of Young diagrams of size $i$. Now the set
\{$\widehat{\widetilde{y}}_\lambda: |\lambda|=n$ or less  than $n$
by an even number\} forms a basis for the subspace
$\widehat{K_n}$, each $\widehat{\widetilde{y}}_\lambda$ is an
eigenvector with a unique eigenvalue, we conclude that the
eigenspaces are all $1$-dimensional.
\end{proof}

This establishes that any element in $\widehat{K_n}$ can be
written  as a linear combination
 of the $\widehat{\widetilde{y}}_\lambda$ with $|\lambda|$ being equal to $n$
 or less than $n$ by an even number. It follows that any element of $\widehat{K_n}$
 which is an eigenvector of $\varphi$ must be a multiple of some
 $\widehat{\widetilde{y}}_\lambda$ as given above.

\section{Kauffman polynomials of some generalized Hopf links}
In this section, we will apply the techniques described in the
previous  sections to show that it is possible to compute the
Kauffman polynomial of some generalized Hopf links. We will give
specific examples in the next section.

Consider $H(n,0)$ in $K(S^1\times D^2)$, then $H(n,0)=A_1^n$
($A_1$  denotes the unknot in the annulus representing the
longitude).

From the previous section, a basis for $\widehat{K_n}$ is given by
the set \{$\widehat{\widetilde{y}}_\lambda: |\lambda|=n$ or less
than $n$ by  an even number\}. We have

$$A_1^n=\sum_{\lambda}d_\lambda \widehat{\widetilde{y}}_\lambda,$$
where the sum is over all Young diagram $\lambda$ with $|\lambda|$
equal to $n$ or less than $n$ by an even number.

Therefore, in $K(S^1\times D^2)$,

$$H(n, m)=\varphi^m(A_1^n)=\sum_{|\lambda|=i}d_\lambda \varphi^m(\widehat{\widetilde{y}}_\lambda)=\sum_{|\lambda|=i}d_\lambda c_{\lambda}^m \widehat{\widetilde{y}}_\lambda.$$

So by evaluating in the plane, the Kauffman polynomial of $H(n, m)$ is:
$$<H(n, m)>=\sum_{|\lambda|=i}d_\lambda c_{\lambda}^m  <\lambda>,$$
where $<\lambda>$ is the Kauffman polynomial of
$\widehat{\widetilde{y}}_\lambda$ given by Wenzl
\cite{W88}[Theorem 5.5] and by Beliakova and Blanchet
\cite{BB2000}[Theorem 7.5].

\section{The Kauffman polynomials of $H(1, n)$ and $H(2, n)$}
\subsection{The Kauffman polynomials of $H(1, n)$}
When $\lambda$ is a Young diagram with a single cell, we have
$c_\lambda=[\delta+(s-s^{-1})(\alpha-\alpha^{-1})]$.

$$<\varphi(\widehat{\widetilde{y}}_\lambda)>=<Hopf link>= \delta [\delta+(s-s^{-1})(\alpha-\alpha^{-1})].$$
Therefore,

$$<H(1,n)>=<\varphi^n(\widehat{\widetilde{y}}_\lambda)>=\delta [\delta+(s-s^{-1})(\alpha-\alpha^{-1})]^n.$$

\subsection{The Kauffman polynomials of $H(2, n)$}
Let $\lambda$ be the Young diagram of size two with the two cells on one row,
let $\mu$ be the Young diagram of size two with the two cells on one column.
$$\lambda=\quad\raisebox{-3mm}{\epsfysize0in\epsffile{2l.ai}}\quad, \quad\quad\mu=\quad\raisebox{-5mm}{\epsfysize0in\epsffile{2b.ai}}.$$
Recall in the Hecke algebra $H_2$, the idempotents corresponding
to  $\lambda$ and $\mu$ are the symmetrizer and antisymmetrizer
denoted by $f_2$ and $g_2$ by Beliakova and Blanchet
\cite{BB2000}[Chapter 2], where
$$f_2=\frac{s^{-1}}{s+s^{-1}}\raisebox{-3mm}{\epsfxsize.3in\epsffile{parra.ai}}\quad + \frac{s^{-1}}{s+s^{-1}}\quad \raisebox{-3mm}{\epsfxsize.3in\epsffile{left.ai}}$$
and
$$g_2=\frac{s}{s+s^{-1}}\raisebox{-3mm}{\epsfxsize.3in\epsffile{parra.ai}}\quad - \frac{s^{-1}}{s+s^{-1}}\quad \raisebox{-3mm}{\epsfxsize.3in\epsffile{left.ai}}.$$

From the definition of $f_2$ and $g_2$, in $H_2$, we have,
$$\raisebox{-3mm}{\epsfxsize.3in\epsffile{parra.ai}}\quad=f_2+g_2.$$

Now apply the multiplicative homomorphism $s_2: H_2\to K_2$,
from the construction of $s_2$ in the proof of
\cite{BB2000}[Theorem 3.1],

$$s_2(\quad \raisebox{-3mm}{\epsfxsize.3in\epsffile{parra.ai}}\quad)=\quad\raisebox{-3mm}{\epsfxsize.3in\epsffile{parra1.ai}}-\delta^{-1} \quad\raisebox{-3mm}{\epsfxsize.3in\epsffile{parra2.ai}}\quad.$$
Therefore,
\begin{alignat}{2}
\raisebox{-3mm}{\epsfxsize.3in\epsffile{parra1.ai}}&=s_2(\quad \raisebox{-3mm}{\epsfxsize.3in\epsffile{parra.ai}}\quad)+  \delta^{-1} \quad\raisebox{-3mm}{\epsfxsize.3in\epsffile{parra2.ai}}\notag\\
&=s_2(f_2)+s_2(g_2)+ \delta^{-1} \quad\raisebox{-3mm}{\epsfxsize.3in\epsffile{parra2.ai}}.\notag
\end{alignat}

Now by embedding this into the skein of the annulus, we have,

\begin{alignat}{2}
H(2,0)&=\widehat{\widetilde{y}_{\lambda}}+\widehat{\widetilde{y}_{\mu}}+\delta^{-1}\delta\phi\notag\\
&=\widehat{\widetilde{y}_{\lambda}}+\widehat{\widetilde{y}_{\mu}}+\phi.\notag
\end{alignat}

So the Kauffman polynomial of $H(2, n)$ is
\begin{alignat}{2}
<H(2,n)>&=<\varphi^n(\widehat{\widetilde{y}_{\lambda}}+\widehat{\widetilde{y}_{\mu}}+\phi)>\notag\\
         &=<\varphi^n(\widehat{\widetilde{y}_{\lambda}})>+<\varphi^n(\widehat{\widetilde{y}_{\mu}})>+<\varphi^n(\phi)>
         \notag \\
 &=c_\lambda^n<\lambda>+c_\mu^n<\mu>+\delta^n ,\notag
\end{alignat}
where
$c_\lambda=\delta+(s-s^{-1})(\alpha(1+s^2)-\alpha^{-1}(1+s^{-2}))$
and
$c_\mu=\delta+(s-s^{-1})(\alpha(1+s^{-2})-\alpha^{-1}(1+s^{2}))$,
and from \cite{BB2000}[Theorem 7.5],
$$<\lambda>=\frac{(\alpha-\alpha^{-1})(\alpha s-\alpha^{-1}s^{-1}+s^2-s^{-2})}{(s+s^{-1})((s-s^{-1})^2}$$
and
$$<\mu>=\frac{(\alpha-\alpha^{-1})(\alpha s^{-1}-\alpha^{-1}s+s^2-s^{-2})}{(s+s^{-1})((s-s^{-1})^2}.$$
As a special case, when $n=0$,
we have $<H(2,0)>=<\lambda>+<\mu>+1=\delta^2$.

\end{document}